\documentstyle{article}
\begin{document}
\newcommand{\ti}[1]{\mbox{$\tilde{#1}$}}
\newcommand{\wti}[1]{\mbox{$\widetilde{#1}$}}
\newcommand{\T}{\mbox{$\widetilde{T}$}}
\newcommand{\U}{\mbox{$\widetilde{U}$}}
\newcommand{\SS}{\mbox{$\widetilde{S}$}}
\newcommand{\nul}{\mbox{$\tilde{0}$}}
\newcommand{\al}{\mbox{$\alpha$}}
\newcommand{\be}{\mbox{$\beta$}}
\newcommand{\dd}{\mbox{$\delta$}}

\newcommand{\petit}{\scriptstyle}
\newcommand{\la}{\mbox{$\lambda$}}
\newcommand{\e}{\mbox{$\epsilon$}}
\newcommand{\norm}[1]{\mbox{$\|#1\|$}}
\newcommand{\YY}{\mbox{$Y'$}}
\newcommand{\TT}{\mbox{$T'$}}
\newcommand{\tutu}{\mbox{$t'$}}
\newcommand{\vv}{\mbox{$v'$}}
\newcommand{\ZZ}{\mbox{$Z'$}}
\newcommand{\yy}{\mbox{$y'$}}
\newcommand{\zz}{\mbox{$z'$}}
\newcommand{\yiy}{\mbox{${\scriptstyle y}'$}}
\newcommand{\ziz}{\mbox{${\scriptstyle z}'$}}

\font\tenBbb=msbm10  \font\sevenBbb=msbm7  \font\fiveBbb=msbm5
\newfam\Bbbfam
\textfont\Bbbfam=\tenBbb \scriptfont\Bbbfam=\sevenBbb
\scriptscriptfont\Bbbfam=\fiveBbb
\def\Bbb{\fam\Bbbfam\tenBbb}
\def\R{{\Bbb R}}
\def\C{{\Bbb C}}
\def\N{{\Bbb N}}
\def\zed{{\Bbb Z}}

\title{OPERATORS ON SUBSPACES OF HEREDITARILY INDECOMPOSABLE
BANACH SPACES}
\author{V.Ferenczi}
\maketitle
\abstract{ We show that if $X$ is a complex hereditarily indecomposable
space, then every operator from a subspace $Y$ of $X$ to $X$
is of the form $\la I+S$, where $I$ is the inclusion
map and $S$ is strictly singular.}

\section{Introduction}
By {\em space (resp. subspace)}, we shall always mean
infinite dimensional space (resp. subspace). A space
$X$ is {\em hereditarily indecomposable}  
 if no two subspaces of $X$ are in a direct sum.
In the whole article, $X$ stands for
 a hereditarily indecomposable
complex Banach space. It was shown in \cite{GM} that every
operator from $X$ to $X$ is of the form 
$\la I+S$, where $I$ is the identity
map and $S$ is strictly singular.
 We generalize this result by showing
that for every subspace $Y$ of $X$, every operator from
$Y$ to $X$ is of the form $\la I+S$, where $I$ is the inclusion
map and $S$ is strictly singular (this was proved in
\cite{GM} in a particular case).

 In fact, it is a consequence
of the following lemma that it is enough to
prove this for every subspace $Y$ with a basis with constant $2$.

\section{Lemma}
\paragraph{Lemma 1}
{\em Let $Y$ be H.I. and $S$ an operator
 from $Y$ to some Banach space. Let
$Z \subset Y$. Then $S$ is strictly singular if and only if
$S_{/Z}$ is strictly singular.}

\underline{Proof}
The direct implication is clear. Now suppose that 
$S_{/Z}$ is strictly singular while $S$ is not. Then there
is a subspace $\YY$ of $Y$, and $c > 0$ such that
$\forall y \in \YY, \norm{S(y)} \geq c \norm{y}$; furthermore,
given $\e < c/(1+\norm{S})$, there is a subspace $\ZZ$ of $Z$ with 
$\norm{S_{/\ZZ}} \leq \e$. As $Y$ is H.I., there exist two unit
vectors
$y \in \YY$, $z \in \ZZ$ with $\norm{y-z} \leq \e$. Then
$ c-\e \leq \norm{S(y-z)} \leq \norm{S}\e $,
a contradiction.

\paragraph{Consequence}
In particular, this lemma can be applied to an operator
from $Y \subset X$ to $X$.
Now suppose that the claim of the article is true for every subspace with a basis
with constant $2$. Let
$Y \subset X$, and $T:Y \rightarrow X$.
 We know that $Y$ contains a normalized basic
sequence with constant $2$;
 let $Z$ be the subspace associated to it.
By hypothesis, there exists $\la$ such that
$T_{/Z}-\la I_{/Z}$ is strictly singular. By the lemma,
$T-\la I_{/Y}$ is still strictly singular.

\section{A filter on block-subspaces}
 
A normalized basic sequence $(y_i)_{i \in \N}$ in $X$ will be 
denoted by $y$. We denote by $[y]$ 
the closed subspace generated by $y$, by $I_y$ the 
 inclusion map from $[y]$ to $X$, by $y^{k}$ the 
normalized basic sequence $(y_{i+k})_{i \in \N}$.
For an operator $T$ such that 
the restriction of $T$ to $[y]$ is defined, $T_{/y}$ stands
for this restriction.
If $(y_i)_{i \in \N}$ is a basic sequence such that
for every $i$, $y_i \neq 0$, then $|y|$ denotes the normalized
basic sequence $(y_i/\norm{y_i})_{i \in \N}$.
By $\yy \subset y$, we mean that $\yy$ is a normalized
block basic sequence of $y$.
Let $e$ be a normalized basic sequence in $X$ with constant~$2$. 
Let $J=\{y / y \subset e\}$.  

\paragraph{Definition 1}

Given $y$ and $z$ in $J$, we say that they
 are
{\em similar}, and write 
{\em $y \sim z$} if
  $\Sigma \norm{z_i-y_i} < +\infty$.
\paragraph{Properties}
Similarity is an equivalence relation.
Furthermore, for every $y \sim z$ and every
$\yy \subset y$, there exists $\zz \subset z$ such that
$\yy \sim \zz$. Indeed, for $k \in \N$, let $\zz_k$
have the same coordinates on $(z_i)$ as $\yy_k$ on $(y_i)$: 
the normalization $|\zz|$ of $(\zz_k)_{k \in \N}$
 is similar to $\yy$.

\paragraph{Definition 2}
{\em For $y$ and $z$ in $J$, we say
that $y \leq z$ if 
$\exists \zz \subset z / y \sim \zz $.}

It is a consequence of the properties of similarity
that $\leq$ is a preordering.

It is also a filter. Indeed let $y$ and
 $z$ in $J$;
 using the fact that $X$
is H.I., it is possible to find $\yy \subset y$ and
$\zz \subset z$ with $\yy \sim \zz$. 
We have that $\yy \leq y$ and $\yy \leq z$.

\paragraph{Definition 3}
We define a filter
on the set of block subspaces of $[e]$
by letting
$Y \leq Z$ if there exists $y \leq z$ with $Y=[y]$, $Z=[z]$.

\paragraph{Definition 4}

Let $\cal U$ be an ultrafilter on $J$. For $y \in J$,
 let $B_{y}$
be the Banach space of bounded operators from $[y]$ to $X$,
 and let
$\norm{.}_y$ be the seminorm 
$\lim_{k \rightarrow +\infty}\norm{._{/y^k}}$.

Let $\cal B$ be the quotient space
 of $l_{\infty}((B_{y})_{y \in J})$ by
the kernel of the seminorm
$\lim_{\cal U}\norm{.}_{y}$.

\paragraph{Lemma 2}
{\em Let $y \in J$ and $z \in J$ be similar . The operator
$p_{yz}$
from $[y]$ to $[z]$ defined by $\forall i \in \N, p_{yz}(y_i)=z_i$
is bounded and 
$\norm{p_{yz}-I_y}_y=0 $.}

\section{A morphism from ${\cal L}(Y,X)$ to $\cal B$}
\paragraph{Definition 5}
Let $T:Y \rightarrow X$ with $Y=[y]$.
 For $z \leq y$, 
 let $\yy \subset y$ such that $z \sim \yy$,
let $T_z$ be the element $T p_{z\yiy}$ of $B_z$.
Let $\T$ be the element of $\cal B$ associated
 to $(T_z)_{z \leq y}$ (the value of $T_z$ for the 
other values of $z$ has no effect on the value of $\T$,
 take for example $T_z=0$).

 This defini\-tion does not depend
 on the choice of $\yy$. Indeed, let $\yy_1~\sim~\yy_2$ be
two choices of $\yy$;
let $T_{iz}$ be asso\-ciated to $\yy_i$; 
then $T_{1z}-T_{2z}=T(p_{z{\petit y}'_1}-p_{z{\petit y}'_2})$ so that
by Lemma~2,
 $\norm{T_{1z}-T_{2z}}_z=0$.

It is easy to check that $z \rightarrow \norm{T_z}_{z}$
 is increasing, so that $\norm{\T}$ is a simple limit.

Let $B$ be the set of elements of $\cal B$ of the form $\T$.

\paragraph{Remark 1}
\

If $T:Y \rightarrow X$ and $U:Y \rightarrow X$, then 
$\wti{T+U}=\T+\U$.

By Lemma 2, for all $y \sim z$, $\ti{p_{yz}}=\ti{I_y}$.

\section{Lemmas} 
\paragraph {Lemma 3}
{\em Let $T:Y \rightarrow X$ with $Y=[y]$, and let $z \leq y$.
Then $\ti{T_z}=\T$.}

\underline{Proof}
Let $\yy \subset y$ such that $z \sim \yy$.
Let $\zz \subset z$.
For $k$ in $\N$, let $t_k=p_{z\yiy}(\zz_k)$.
It has the same coordinates on $(\yy_i)$ as 
$\zz_k$ on $(z_i)$, so $\zz \sim |t|$.
We have

\[ \sum_k |\norm{t_k}-1| \leq \sum_k \norm{t_k-\zz_k}
\leq 4 \sum_i \norm{\yy_i-z_i}, \]

so $\sum_k |\norm{t_k}-1| $ converges.

Furthermore, defining $T_{\ziz}$ as $Tp_{\ziz |t|}$ and 
$T_z$ as $Tp_{z\yiy}$, we get for every~$k$:

\[ \norm{(T_{\ziz}-(T_z)_{\ziz})(\zz_k)} =
    \norm{T(|t|_k-t_k)} \leq
\norm{T} | \norm{t_k}-1 | .\]

These two points imply that
 $\norm{T_{\ziz}-(T_z)_{\ziz}}_{\ziz}=0$;
 so $\ti{T_z}=\T$.

\paragraph{Consequences}
\

If $\yy \subset y$ then $\wti{T_{/\yiy}}=\T$.

Let $\T$ and $\U$ be in $B$, with $T:Y \rightarrow X$
and $U:Z \rightarrow X$; let $t \leq y$ and $t \leq z$;
then $\T+\U=\T_t+\U_t=\wti{T_t+U_t}$. This proves
that $B$ is a linear space.

\

>From now on, $Y$ (resp. $\YY$, $Z$, $\ZZ$)
 stands for a block subspace $[y]$ 
(resp. $[\yy]$, $[z]$, $[\zz]$).
 Recall that $Y \leq Z$ means $y \leq z$.

\paragraph{Proposition 1}
{\em Let $S:Y \rightarrow X$. Then the following properties
are equivalent:

(a) $\SS=\nul$.

(b) $S$ is strictly singular.}

\underline{Proof}

(a) means
$ \forall \e > 0, \exists z \leq y
 \ s.t.\ \norm{S_{z}}_{z} \leq \e$
or equivalently, the assertion:
 \[(c)\ \forall \e > 0, \exists \yy \subset y,
 \exists k\ s.t.\ \norm{S_{\yiy^k}} \leq 2\e .\]
This last assertion implies (b) and (b) implies (c) is shown to be true
in \cite{LT}.

\paragraph {Lemma 4} {\em Let $T:Y \rightarrow X$.
 Let $Z$ be a block subspace of $X$.
There exists $\YY \subset Y$, $\TT: \YY \rightarrow Z$ with
$\T=\wti{\TT}$.}

\underline{Proof} First notice that if $S:Y \rightarrow X$ is such
 that $\Sigma \norm{S(y_i)}$ converges, then $\SS=\nul$.

Now let $T:Y \rightarrow X$. If $\T=\nul$ then one 
can choose $\TT=0$, $\YY=Y$. Suppose now $\T \neq \nul$.
By the first consequence of Lemma~3 and Proposition~1, passing
to a subspace,
 we may assume that $T$ is an isomorphism.

By projections on big enough intervals,
selecting
 block vectors, and using repeatedly the fact that $X$ is H.I.,
it is possible to build a normalized block
 basic sequence $\yy \subset y$ and a (non-normalized)
block basic sequence $(z_i)_{i \in \N}$ in $Z$ such that
the sequence $\Sigma \norm{T(\yy_i)-z_i}$ converges.

The operator $\TT$ defined on $[\yy]$ by
$\TT(\yy_i)=z_i$ satisfies $\ti{\TT}=\T$ by the first remark,
and takes its values in $Z$.

\section{A product on $B$}
\paragraph{Definition 6}
Let $\T \in B$ with $T:Z \rightarrow X$.
Let $\U \in B$. Let $V \leq Z$. By Lemma 4,
 we can assume that $\U$
is associated to $U:Y \rightarrow V$.
We want to define $\T\U$ as $\wti{T_v U}$. We 
need to show that $\wti{T_v U}$ does not depend on the choice of 
$T$ and $U$.

First, for $i=1,2$, let $V_i \leq Z$, and
 $U_i:Y_i \rightarrow V_i$ be such that
$\ti{U_i}=\U$. Let $\zz_i \subset z$ with $v_i \sim \zz_i$.
Let $\yy \leq y_1$ and $\yy \leq y_2$.
Then \[\wti{T_{v_2}U_2}-\wti{T_{v_1}U_1}=
     \wti{(T_{v_2}U_2)_{\yiy}}-\wti{(T_{v_1}U_1)_{\yiy}}=
     \wti{T_{v_2}U_{2\yiy}-T_{v_1}U_{1\yiy}}.\]

Now \[ T_{v_2}U_{2\yiy}-T_{v_1}U_{1\yiy}=
      T(p_{v_2 {{\petit z}'}_2} U_{2\yiy}-p_{v_1 {{\petit z}'}_1} U_{1\yiy})\]

\[  T_{v_2}U_{2\yiy}-T_{v_1}U_{1\yiy}     =
 T[(p_{v_2 {{\petit z}'}_2}-I_{v_2})U_{2\yiy}
+(U_{2\yiy}-U_{1\yiy})+(I_{v_1}-p_{v_1 {{\petit z}'}_1}) U_{1\yiy}].\]

Using Remark 1, and the fact that the space of strictly singular
operators is a two-sided ideal, we get that
\[\wti{T_{v_2}U_2}-\wti{T_{v_1}U_1}=\nul.\]

We now prove that $\wti{T_v U}$ does not depend on the choice of
$T$.
Let for $i=1,2$ $T_i:Z_i \rightarrow X$ be such that
$\ti{T_i}=\T$. Let $V \leq Z_1$ and $V \leq Z_2$ and
$U:Y \rightarrow V$ be a representative for $\U$.
Then $\wti{(T_2)_v U}-\wti{(T_1)_v U}=\wti{((T_2)_v-(T_1)_v)U}=\nul$.

So we can define $\T \U$ as $\wti{T_v U}$ without ambiguity.

\paragraph{Remark 2}
If $T:Z \rightarrow X$ and $U:Y \rightarrow Z$, then
$\T \U= \wti{TU}$.
Using appropriate representatives one can then
show that $B$ is an algebra; in particular $\ti{1}$ the common
value associated to all $I_y$ for $y \in J$ is neutral for
the multiplication. 

\paragraph{Proposition 2}
{\em $B$ is a field.}

\underline{Proof}
Let $\T \neq \nul$, and let us choose a representative
 $T:Y \rightarrow Z$ such that $Z$ is a block subspace.
By Proposition 1 and the first consequence of Lemma 3,
 we may assume that $T:Y \rightarrow Z$
is an isomorphism onto.
Then $\T \wti{T^{-1}} = \wti{T T^{-1}}= \wti{I_{/Z}}=\ti{1}$; in the
same way $\wti{T^{-1}} \T = \ti{1}$; so $\T$ is invertible.
As $B$ is not a singleton ($\ti{1} \neq \nul$), $B$ is a field.

\section{$B$ is a Banach algebra.}

\paragraph{Lemma 5}

{\em For all $\la \in B$,}
 \[ \norm{\la}=
\inf\{\norm{T}: y \in J,\ T \in {\cal L}(Y,X),\ \T=\la\}.\]

\underline{Proof} If $\la= \nul$, the result is clear. Now
assume $\la \neq \nul$. For every $Y$, if $T$ belongs to ${\cal L}(Y,X)$ then
$\norm{\T} \leq \norm{T_y}_y \leq \norm{T}$,
 so that
$\norm{\la} \leq \inf_{\T=\la}\norm{T}$.

Furthermore, let $T:Z \rightarrow X$
 satisfy $\la=\T$. For every $\e > 0$, there exists $\zz \subset z$
so that $\norm{T_{\ziz}}_{\ziz} \leq (1+\e)\norm{\la}$.
 So there exists $k$ 
such that $\norm{T_{\ziz^k}} \leq (1+\e)^2 \norm{\la}$.
As $\wti{T_{\ziz^k}}=\T$, we have $ \inf_{\T=\la} \norm{T} \leq  \norm{\la}$.

\underline{Remark}
We say that a representative $T$ for $\T$ is {\em $\e$-minimal} if 
$\norm{T} \leq (1+\e)\norm{\T}$.
 We have proved that if $T$ belongs to ${\cal L}(Y,X)$, then
for any $Z \leq Y$ and any $\e > 0$,
we can find an $\e$-minimal representative for $\T$ defined on a block subspace
of $Z$.

\paragraph{Proposition 3}
{\em The norm on $B$ is an algebra norm.}

\underline{Proof}
 Let $\T$ and $\U$ be in $B$.
Let $\e > 0$.
Let $T:Y \rightarrow X$ be an $\e$-minimal representative 
for $\T$. Choose an $\e$-minimal representative for $\U$ from
some $Z$ into $Y$ (it is possible by the previous remark and
 Lemma 4).

Then $\norm{\T\U}=\norm{\wti{TU}}
 \leq \norm{TU} \leq \norm{T}\norm{U} 
\leq (1+\e)^2 \norm{\T}\norm{\U}$.

As $\e$ is arbitrary, the norm on $B$ is an algebra norm.

\paragraph{Proposition 4}
 {\em $B$ is a Banach space.}

\underline{Proof}
We show that any normally conver\-ging series conver\-ges in $B$.

Let $(\la_n)_{n \in \N}$ be a normally converging series.
Let $T_0:Y_0 \rightarrow X$ be a $1$-minimal represen\-tative
 for $\la_0$.
Given $T_{n-1}:Y_{n-1} \rightarrow X$, by the remark at the
end of Lemma~5, we can find a $1$-minimal
representative $T_n:Y_n \rightarrow X$ for $\la_n$ with 
$Y_n$ a block subspace of $Y_{n-1}$.

Let $y \in J$ such that $\forall n, y_n \in Y_n$.
We define an operator $H$ on $Y$ by
 $H(y_n)=\sum_{k=0}^n T_k(y_n)$.
Let us evaluate $\norm{(H-\sum_{k=0}^n T_k)_{/y^n}}$.

Let $u=\sum_{i=n}^{+\infty}a_i y_i$ be a vector
 in $[y^n]$. Then

\[ \left(H-\sum_{k=0}^n T_k \right)(u)=
\sum_{i=n}^{+\infty}
 a_i \left(\sum_{k=0}^{i}T_k(y_i)-\sum_{k=0}^{n} T_k(y_i) \right). \]

\[ \left(H-\sum_{k=0}^n T_k\right)(u)=
\sum_{i=n}^{+\infty}
a_i \left(\sum_{k=n+1}^{i}T_k(y_i) \right). \]

\[\left(H-\sum_{k=0}^n T_k \right)(u)=
\sum_{k=n+1}^{+\infty} T_k 
\left(\sum_{i=k}^{+\infty}a_i y_i \right). \]

So \[\norm{(H-\sum_{k=0}^n T_k)(u)}
 \leq 3 \left(\sum_{k=n+1}^{+\infty} \norm{T_k}\right) \norm{u}.\]

In parti\-cular, $H$ is a bounded operator.
 Furthermore,
\[ \lim_{n \rightarrow +\infty}
 \norm{(H-\sum_{k=0}^n T_k)_{/y^n}}=0.\]
By Lemma 5,
 $\lim_{n \rightarrow +\infty}
\norm{\ti{H}-\sum_{k=0}^{n}\la_k}=0$. This proves that
$\sum_{k=0}^{+\infty} \la_k$ converges in $B$ to $\ti{H}$.

\section{Conclusion}

The algebra $B$ is a Banach algebra and also a field. 
By Gel'fand-Mazur theorem, we can identify it to the field of complex numbers.

Now let $T$ belong to ${\cal L}(Y,X)$.
Then $\wti{T-\T I_{/Y}}=\nul$, that is  $T-\T I_{/Y}$
is strictly singular.

So for any block subspace $Y$ of $[e]$, all operators on $Y$ are of
the required form. As $e$ is arbitrary, this is true for any
subspace with a basis with constant~$2$, and thus for any subspace.

\

\

\

\

Equipe d'Analyse et de Math\'{e}matiques Appliqu\'{e}es,

Universit\'{e} de Marne-la-Vall\'{e}e,

2, rue de la Butte Verte,

93166 Noisy le Grand Cedex, France.

ferenczi@math.univ-mlv.fr.

\end{document}